\documentclass[onecolumn]{article}

\usepackage{graphicx,epsfig}
\usepackage{times}
\usepackage{hyperref}
\usepackage{amssymb, amsthm}
\usepackage {amssymb,latexsym,epic,eepic,stmaryrd,url}
\newtheorem{thm}{Theorem}
\newtheorem{clm}{Claim}

\newtheorem{fact}{Fact}
\newtheorem{lem}{Lemma}
\newtheorem{asm}{Assumption}

\newtheorem{cor}{Corollary}

\begin{document}
\title{\Large \bf Acyclic Edge Coloring of Graphs with maximum degree 4}

\author{Manu Basavaraju\thanks{Computer Science and Automation department,
Indian Institute of Science,
Bangalore- 560012,
India.  {\tt manu@csa.iisc.ernet.in}} \and L. Sunil Chandran\thanks{ {\bf Corresponding Author:} Computer Science and Automation department,
Indian Institute of Science,
Bangalore- 560012,
India.  {\tt sunil@csa.iisc.ernet.in}}
}

\date{}
\pagestyle{plain}
\maketitle

\begin{abstract}

An $acyclic$ edge coloring of a graph is a proper edge coloring such that there are no bichromatic cycles. The \emph{acyclic chromatic index} of a graph is the minimum number k such that there is an acyclic edge coloring using k colors and is denoted by $a'(G)$. It was conjectured by Alon, Sudakov and Zaks that for any simple and finite graph $G$, $a'(G)\le \Delta+2$, where $\Delta =\Delta(G)$ denotes the maximum degree of $G$. We prove the conjecture for connected graphs with $\Delta(G) \le 4$, with the additional restriction that $m \le 2n-1$, where $n$ is the number of vertices and $m$ is the number of edges in $G$. Note that for any graph $G$, $m \le 2n$, when $\Delta(G) \le 4$. It follows that for any graph $G$ if $\Delta(G) \le 4$, then $a'(G) \le 7$.

\end{abstract}


\section{Introduction}

All graphs considered in this paper are finite and simple. A proper \emph{edge coloring} of $G=(V,E)$ is a map $c: E\rightarrow S$ (where $S$ is the set of available $colors$ ) with $c(e) \neq c(f)$ for any adjacent edges $e$,$f$. The minimum number of colors needed to properly color the edges of $G$, is called the chromatic index of $G$ and is denoted by $\chi'(G)$. A proper edge coloring c is called acyclic if there are no bichromatic cycles in the graph. In other words an edge coloring is acyclic if the union of any two color classes induces a set of paths (i.e., linear forest) in $G$. The \emph{acyclic edge chromatic number} (also called \emph{acyclic chromatic index}), denoted by $a'(G)$, is the minimum number of colors required to acyclically edge color $G$. The concept of \emph{acyclic coloring} of a graph was introduced by Gr\"unbaum \cite{Grun}. The \emph{acyclic chromatic index} and its vertex analogue can be used to bound other parameters like \emph{oriented chromatic number} and \emph{star chromatic number} of a graph, both of which have many practical applications, for example, in wavelength routing in optical networks ( \cite{ART}, \cite{KSZ} ). Let $\Delta=\Delta(G)$ denote the maximum degree of a vertex in graph $G$. By Vizing's theorem, we have $\Delta \le \chi'(G) \le \Delta +1 $(see \cite{Diest} for proof). Since any acyclic edge coloring is also proper, we have $a'(G)\ge\chi'(G)\ge\Delta$. \newline

It has been conjectured by Alon, Sudakov and Zaks \cite{ASZ} that $a'(G)\le\Delta+2$ for any $G$. Using probabilistic arguments Alon, McDiarmid and Reed \cite{AMR} proved that $a'(G)\le60\Delta$. The best known result up to now for arbitrary graph, is by Molloy and Reed  \cite{MolReed} who showed that $a'(G)\le16\Delta$. Muthu, Narayanan and Subramanian \cite{MNS1} proved that $a'(G)\le4.52\Delta$ for graphs $G$ of girth at least 220 (\emph{Girth} is the length of a shortest cycle in a graph).\newline

Though the best known upper bound for general case is far from the conjectured $\Delta+2$, the conjecture has been shown to be true for some special classes of graphs. Alon, Sudakov and Zaks \cite{ASZ} proved that there exists a constant $k$ such that $a'(G)\le\Delta+2$ for any graph $G$ whose girth is at least $k\Delta\log\Delta$. They also proved that $a'(G)\le\Delta+2$ for almost all $\Delta$-regular graphs. This result was improved by Ne\v set\v ril and Wormald \cite{NesWorm} who showed that for a random $\Delta$-regular graph $a'(G)\le \Delta+1$. Muthu, Narayanan and Subramanian proved the conjecture for grid-like graphs \cite{MNS2} and outer planar graphs \cite{MNS3}. In fact they gave a better bound of $\Delta+1$ for those classes of graphs. From Burnstein's \cite{Burn} result it follows that the conjecture is true for subcubic graphs. Skulrattankulchai \cite{Skul} gave a polynomial time algorithm to color a subcubic graph using $\Delta+2 = 5$ colors. Recently Basavaraju and Chandran \cite{MBSC} proved that connected non-regular subcubic graph can be acyclically edge colored using $4$ colors.

Determining $a'(G)$ is a hard problem both from a theoretical and from an algorithmic point of view. Even for the simple and highly structured class of complete graphs, the value of $a'(G)$ is still not determined exactly. It has also been shown by Alon and Zaks \cite{AZ} that determining whether $a'(G)\le3$ is NP-complete for an arbitrary graph $G$. The vertex version of this problem has also been extensively studied ( see \cite{Grun}, \cite{Burn}, \cite{Boro}). A generalization of the acyclic edge chromatic number has been studied. The \emph{r-acyclic edge chromatic number} $a'_r(G)$ is the minimum number of colors required to color the edges of the graph $G$ such that every cycle $C$ of $G$ has at least min\{$\vert C \vert$,$r$\} colors ( see \cite{GeRa}, \cite{GrePi}).

~~~~~~~~~

\noindent\textbf{Our Result:}In this paper, we prove the following Theorem:

\begin{thm}
\label{thm:thm1}
Let $G$ be a connected graph on $n$ vertices, $m \le 2n-1$ edges and maximum degree $\Delta \le 4$, then $a'(G)\le 6$. $($Note that if $\Delta(G) \le 4$, then $m \le 2n$ always$)$.
\end{thm}

Note that if $\Delta(G) < 4$, it follows from Burnstein's \cite{Burn} result that $a'(G) \le 5$ and is tight.

\begin{cor}
\label{cor:cor1}
Let $G$ be a graph with maximum degree $\Delta \le 4$. Then $a'(G) \le 7$.
\end{cor}
\begin{proof}
If $\Delta(G) \le 4$, then $m \le 2n$ for each connected component. If $m \le 2n-1$, by Theorem \ref{thm:thm1} $a'(G) \le 6$ for each connected component. Otherwise if $m=2n$, we can remove an edge from each connected component and color the resulting graph with at most 6 colors. Now the removed edges of each component could be colored using a new color. Thus $a'(G) \le 7$.
\end{proof}

\noindent \textbf{Remark:} There exists graphs with $\Delta(G) \le 4$ that require at least 5 colors to be acyclically edge colored. For example, any graph with $\Delta(G) = 4$ and $m=2n-1$ require 5 colors. But we do not know whether there exist any graph with $\Delta(G) \le 4$ that needs 7 colors or even 6 colors to be acyclically edge colored. Thus we feel that the bound of $Corollary$ \ref{cor:cor1} and $Theorem$ \ref{thm:thm1} can be improved. Our proof is constructive and yields an efficient polynomial time algorithm.

\section{Preliminaries}

Let $G=(V,E)$ be a simple, finite and connected graph with maximum degree $4$. Let $x \in V$. Then $N_{G}(x)$ will denote the neighbours of $x$ in $G$. For an edge $e \in E$, $G-e$ will denote the graph obtained by deletion of the edge $e$. For $x,y \in V$, when $e=(x,y)=xy$, we may use $G-\{xy\}$ instead of $G-e$. Let $c:E\rightarrow \{1,2,\ldots,k\}$ be an \emph{acyclic edge coloring} of $G$. For an edge $e\in E$, $c(e)$ will denote the color given to $e$ with respect to the coloring $c$. For $x,y \in V$, when $e=(x,y)=xy$ we may use $c(x,y)$ instead of $c(e)$.

To prove the main result, we plan to use induction on the number of edges. Let $G =(V,E)$ be a graph on $m$ edges where $m \ge 1$. We will remove an edge from $G$ and get a graph $G'=(V,E')$ on smaller number of edges. By induction hypothesis $G'$ will have an acyclic edge coloring $c:E'\rightarrow \{1,2,\ldots,6\}$. Throughout the paper we will consistently assume that the edge we remove from $G$ to get $G'$ is $xy=(x,y)$. Then $E' = E - \{xy\}$. Our intention will be to extend the acyclic edge coloring $c$ of $G'$ to $G$ by assigning an appropriate color for the edge $xy$.

An ($\alpha$,$\beta$)-maximal bichromatic path with respect to an acyclic edge coloring $c$ of $G'$ is a path consisting of edges that are colored using the colors $\alpha$ and $\beta$ alternatingly. An ($\alpha$,$\beta$,$a$,$b$)-maximal bichromatic path is an ($\alpha$,$\beta$)-maximal bichromatic path which starts at the vertex $a$ with an edge colored $\alpha$ and ends at $b$. We emphasize that the edge of the ($\alpha$,$\beta$,$a$,$b$)-maximal bichromatic path incident on vertex $a$ is colored $\alpha$ and the edge incident on vertex $b$ can be colored either $\alpha$ or $\beta$. Thus the notations ($\alpha$,$\beta$,$a$,$b$) and ($\alpha$,$\beta$,$b$,$a$) have different meanings. The following fact is obvious from the definition of acyclic edge coloring:

\begin{fact}
\label{fact:fact1}
Given a pair of colors $\alpha$ and $\beta$ of an acyclic edge coloring $c$, there can be at most one maximal ($\alpha$,$\beta$)-bichromatic path containing a particular vertex $v$, with respect to $c$.
\end{fact}

We denote the set of colors in the acyclic edge coloring $c$ by $C = \{1,2,\ldots,6\}$. For any vertex $u \in V(G')$, we define $F_u =\{c(u,z) \vert z \in N_{G'}(u)\}$. For an edge $ab \in E'$, we define $S_{ab} = F_b - \{c(a,b)\}$. Note that $S_{ab}$ need not be the same as $S_{ba}$.

A color $\alpha$ is a \emph{candidate} for an edge \emph{e} in $G$ with respect to a coloring of $G-e$ if none of the adjacent edges of \emph{e} are colored $\alpha$. A candidate color $\alpha$ is \emph{valid} for an edge \emph{e} if assigning the color $\alpha$ to \emph{e} does not result in any bichromatic cycle in $G$.

Let $F=F_x \cup F_y$. Recall that the vertices $x$ and $y$ are non-adjacent in $G'$. Note that any color $\beta \in C-F$ is a candidate color for the edge $xy$ in $G$ with respect to the coloring $c$ of $G'$. But $\beta$ may not be valid. What may be the reason? It is clear that color $\beta$ is not $valid$ if and only if there exists $\alpha \neq \beta$ such that a ($\alpha$,$\beta$)-bichromatic cycle gets formed if we assign color $\beta$ to the edge $xy$. In other words, if and only if, in $G'$ there existed a ($\alpha$,$\beta$,$x$,$y$) maximal bichromatic path with $\alpha$ being the color given to the first and last edge of this path. Such paths play an important role in our proof. We call them $critical\ paths$. It is formally defined below:

~~~~~~~~~

\noindent\textbf{Critical Path:} For an edge $ab \in E$ an $(\alpha,\beta,$a$,$b$)$ maximal bichromatic path which starts out from a vertex $a$ via an edge colored  $\alpha$ and ends at vertex $b$ via an edge colored $\alpha$ is called an $(\alpha,\beta,ab)$ critical path.

\begin{lem}
\label{lem:lem1}
A candidate color for an edge $e=uv$, is valid if $(F_u \cap F_v) - \{c(u,v)\} = (S_{uv} \cap S_{vu}) = \emptyset$.
\end{lem}
\begin{proof}
Any cycle containing the edge $uv$ will also contain an edge incident on $u$ (other than $uv$) as well as an edge incident on $v$ (other than $uv$). Clearly these two edges are colored differently since $(S_{uv} \cap S_{vu}) = \emptyset$. Thus the cycle will have at least 3 colors and therefore any of the candidate color for the edge $uv$ is valid.
\end{proof}

An obvious strategy to extend the coloring $c$ of $G'$ to $G$ would be to try to assign one of the candidate colors in $C-F$ to the edge $xy$. The condition that a candidate color is not valid for the edge $xy$ is captured in the following fact.

\begin{fact}
\label{fact:fact2}
The color $\beta \in C-F$ is not a $valid$ color for the edge $xy$ if and only if $\exists \alpha \in F_x \cap F_y$ such that there is a $(\alpha,\beta,xy)$  critical path in $G'$.
\end{fact}

If none of the colors in $C-F$ is valid for the edge $xy$, then we can group the colors in $C-F$ into two categories namely $weak$ and $strong$.

~~~~~~~~~

\noindent\textbf{Weak Color:} A color $\beta \in C-F$ is called $weak$ if it forms only one critical path with $x$ and $y$ as end points. Equivalently, there exists only one $\alpha \in F_x \cap F_y$ such that there is a $(\alpha,\beta,xy)$ critical path. A weak color $\beta$ is said to be \emph{actively present} in a set $S_{xa}$, if $\exists k \in N_{G'}(a)$, $k \neq x$ such that $c(a,k)= \beta$ and $(\alpha,\beta,xy)$  critical path contains the edge $(a,k)$. If a weak color $\beta \in S_{xa}$ is not actively present in $S_{xa}$ then it is said to be \emph{passively present} in $S_{xa}$.

\noindent\textbf{Strong Color:}If the color $\beta \in C-F$ is not $weak$, it is called $strong$.\newline

If there are weak colors, it makes sense to try to break the critical path containing one of the weak colors, thus enabling us to use that weak color for the edge $xy$. For this purpose we introduce the concept of $Recoloring$.

~~~~~~~~~

\noindent\textbf{Recolor:} We define $c' = Recolor(c,e,\gamma)$ as the recoloring of the edge $e$ with a candidate color $\gamma$ to get a modified coloring $c'$ from $c$, i.e., $c'(e)= \gamma$ and $c'(f)= c(f)$, for all other edges $f$ in $G'$. The recoloring is said to be proper, if the coloring $c'$ is proper. The recoloring is said to be $acyclic$ ($valid$), if in coloring $c'$ there exists no bichromatic cycle.

Recall that our strategy is to extend the coloring of $G'$ to $G$ by assigning a valid color for the edge $xy$. When all the candidate colors of $xy$ turn out to be $invalid$, we try to $slightly\ modify$ the coloring $c$ of $G'$ in such a way that with respect to the modified coloring, we have a valid color for edge $xy$. \emph{Recoloring} of an edge in the critical path which contained a weak color is one such strategy. Sometimes we resort to a slightly more sophisticated strategy to modify the coloring namely $color\ exchange$ defined below.

~~~~~~~~

\noindent \textbf{Color Exchange:} Let $u,i,j \in V(G')$ and $ui,uj \in E(G')$. We define $c'= ColorExchange(c,ui,uj)$ as the the modification of the current coloring $c$ by exchanging the colors of the edges $ui$ and $uj$, i.e., $c'(u,i)=c(u,i)$, $c'(u,i)=c(u,i)$ and $c'(e)=c(e)$ for all other edges $e$ in $G'$. The color exchange with respect to the edges $ui$ and $uj$ is said to be proper if the coloring obtained after the exchange is proper. The color exchange with respect to the edges $ui$ and $uj$ is $valid$ if and only if the coloring obtained after the exchange is acyclic.

In our proof we use the strategy of color exchange many times and in different contexts. All these contexts are more or less similar but differ in minor details. We would like to capture all these different contexts in a general framework. The configuration defined below is an attempt to formalize this:

~~~~~~~~

\noindent \textbf{Configuration A}
Let $u$ be a vertex and $i,j \in N_{G'}(u)$. Let $N'_{G'}(u) \cup N''_{G'}(u)$ be a partition of $N_{G'}(u)-\{i,j\}$, i.e.,$N'_{G'}(u) \cup N''_{G'}(u)$ = $N_{G'}(u)-\{i,j\}$ and $N'_{G'}(u) \cap N''_{G'}(u) = \emptyset$. The 5-tuple $(u,i,j,N'_{G'}(u),N''_{G'}(u))$ is in $configuration\ A$ if $(u,i),(u,j) \in E'$ and
\begin{enumerate}
 \item $c(u,i) \notin S_{uj}$ and $c(u,j) \notin S_{ui}$
\item $\forall z \in N'_{G'}(u)$, $c(u,z) \notin S_{ui}$ and $c(u,z) \notin S_{uj}$
\end{enumerate}

Suppose $(u,i,j,N'_{G'}(u),N''_{G'}(u))$ is in $configuration\ A$ with respect to the coloring $c$. Let $c'$ be the coloring obtained after the color exchange with respect to the edges $ui$ and $uj$. Then note that condition 1 guarantees that the color $c(u,i)$ is a candidate for edge $uj$ and the color $c(u,j)$ is a candidate for edge $ui$ and thus the coloring obtained after the color exchange is proper. Condition 2 inhibits the possibility of any $(c(u,i),c(u,z))$ or $(c(u,j),c(u,z))$ $bichromatic\ cycles$ being formed for any $z \in N'_{G'}(u)$. Its obvious that there can not be any $(c(u,j),c(u,i))$ $bichromatic\ cycles$ after exchange. Thus the following fact is easy to verify:

\begin{fact}
\label{fact:fact3}
Let the 5-tuple $(u,i,j,N'_{G'}(u),N''_{G'}(u))$ be in $configuration\ A$. Then $c'= ColorExchange(c,ui,uj)$ is not $valid$ if and only if $\exists h \in N''_{G'}(u)$ such that after the color exchange (i.e., in $c'$) there exists an $(\alpha, \beta)$ bichromatic cycle that passes through $h$ for $\alpha \in \{c'(u,i),c'(u,j)\}$ and $\beta=c'(u,h)$.
\end{fact}

\noindent In view of $Fact$ \ref{fact:fact3}, the following $Fact$ is obvious:

\begin{fact}
\label{fact:fact4}
Let the 5-tuple $(u,i,j,N'_{G'}(u),N''_{G'}(u))$ be in $configuration\ A$. Then if $N''_{G'}(u)=\emptyset$, the color exchange $c'= ColorExchange(c,ui,uj)$ is $valid$.
\end{fact}

\begin{lem}
\label{lem:lem2}
Let $u,i,j,a,b \in V(G)$, $ui,uj \in E'$ and $ab \in E$. Also let $\{\alpha,\beta\} \cap \{c(u,i),c(u,j)\} \neq \emptyset$ and $\{i,j\} \cap \{a,b\} = \emptyset$. Suppose there exists an ($\alpha$,$\beta$,$ab$)-critical path that passes through vertex $u$, with respect to the coloring $c$ of $G'$. Let $c'= ColorExchange(c,ui,uj)$ be proper. Then with respect to the coloring $c'$, there will not be any ($\alpha$,$\beta$,$ab$)-critical path in G'.
\end{lem}
\begin{proof}
Note that since we are assuming that the color exchange is proper, $c(u,j) \notin S_{ui}$. Thus $\{\alpha,\beta\} \neq \{c(u,i),c(u,j)\}$ because any $(c(u,i),c(u,j),ab)$ critical path through vertex $u$ will have to involve the edges $ui$ and $uj$. Since $i \notin \{a,b\}$, color $c(u,j) \in S_{ui}$, a contradiction. Let P be the ($\alpha$,$\beta$,$ab$)-critical path. Without loss of generality assume that $\gamma = c(u,i) \in \{\alpha,\beta\}$. Since vertex $u$ is contained in path $P$, we claim that the edge $ui$ belongs to the path $P$. This is because $\gamma = c(u,i) \in F_u$ and hence path $P$ has to involve edge ui by the maximality of path $P$. Let us assume without loss of generality that path $P$ starts at vertex $a$ and reaches vertex $i$ before it reaches vertex $u$. With respect to the coloring $c'$, there will not be any edge adjacent to vertex $i$ that is colored $\gamma$. So the ($\alpha$,$\beta$) maximal bichromatic path that starts at vertex $a$, should end at vertex $i$. Since $i \neq b$, by $Fact$ \ref{fact:fact1} we infer that the ($\alpha$,$\beta$,$ab$) critical path does not exist.
\end{proof}

\section{proof of Theorem 1}
\begin{proof}
We prove the Theorem by induction on the number of edges. Let $H=(V_{H},E_{H})$ be a connected graph of $n$ vertices and $m \le 2n-1$ edges and $\Delta(H) \le 4$. Let the Theorem be true for all connected graphs $W$ such that $\Delta(W) \le 4$ and $\vert E(W)\vert \le 2\vert V(W)\vert-1$, with at most $m-1$ edges. Without loss of generality we can assume that $H$ is 2-connected, since if there are cut vertices in $H$, the acyclic edge coloring of the blocks $B_1,B_2\ldots B_k$ of $H$ can easily be extended to $H$ (Note that each block satisfies the property that $\Delta(B_i) \le 4$ and $\vert E(B_i)\vert \le 2\vert V(B_i)\vert-1$). Thus $\delta(H)\ge2$. Now since $H$ has at most $2n-1$ edges, there is a vertex $x$ of degree at most 3.

Let $y \in N_{H}(x)$. The degree of $y$ is at most 4. Let $H'=H-\{xy\}$, i.e.,$H'=(V_{H'},E_{H'})$, where $V_{H'}=V_{H}$ and $E_{H'}=E_{H}-\{xy\}$. Thus in $H'$, $degree(x) \le 2$ and $degree(y) \le 3$. Note that since $H$ is 2-connected, $H'$ is connected

To avoid certain technicalities in the presentation of the proof, we construct the graph $G'$ from $H'$ as below. If $degree_{H'}(x)=2$, $degree_{H'}(y)=3$ and $\forall z \in N_{H'}(x)\cup N_{H'}(y), degree_{H'}(z)=4$, then let $G'=H'$ and $G=H$. Otherwise, we construct the graph $G'=(V',E')$ from $H'$ in the following manner. First add pendant vertices as neighbours to the vertices $x$ and $y$ such that $degree_{G'}(x)=2$ and $degree_{G'}(y)=3$. Next add pendant vertices as neighbours to the newly added vertices and $\forall z \in N_{H'}(x)\cup N_{H'}(y)$ such that $\forall z \in N_{G'}(x)\cup N_{G'}(y), degree_{G'}(z)=4$. Note that since $H'$ was connected, $G'$ is also connected. Let $G = G'\cup \{xy\}$, i.e., $G=(V,E)$, where $V=V'$ and $E=E'+\{xy\}$.

By induction hypothesis, graph $H'$ is acyclically edge colorable using $6$ colors. Note that we can easily extend the coloring of $H'$ to $G'$ by coloring each of the newly added edges with the available colors satisfying the acyclic edge coloring property. Let $c_0 : E'\rightarrow \{1,2,.....,6\}$ be a acyclic edge coloring of $G'$. It is easy to see that if we extend the acyclic edge coloring of $G'$ to $G$ by assigning an appropriate color to the edge $xy$, then this coloring also corresponds to the acyclic edge coloring of $H$, since $H$ is a subgraph of $G$.

Our intention will be to extend the acyclic edge coloring $c_0$ of $G'$ to $G=G'+\{xy\}$ by assigning an appropriate color for the edge $xy$. We denote the set of colors of $c_0$ by $C = \{1,2,3,4,5,6\}$.

Let $N_{G'}(x)=\{a,b\}$ and $N_{G'}(y)=\{a',b',d'\}$. Note that $N_{G'}(x) \cap N_{G'}(y)$ need not be empty. Also recall that $degree_{G'}(a)=degree_{G'}(b)=4$. Let $N_{G'}(a)=\{x,k_1,k_2,k_3\}$ and $N_{G'}(b)=\{x,l_1,l_2,l_3\}$.

\subsection*{case 1: $F_x \cap F_y = \emptyset$}
Since $\vert F \vert = 5$, $\vert C-F \vert = 1$. Clearly the $candidate$ color in $C-F$ is valid for the edge $xy$.

\subsection*{case 2: $\vert F_x \cap F_y\vert = 2$}
\begin{asm}
\label{asm:asm1}
Without loss of generality let $F_x = \{1,2\}$ and $F_y = \{1,2,3\}$. Thus $F=\{1,2,3\}$.
\end{asm}

By $Assumption$ \ref{asm:asm1}, $C-F = \{4,5,6\}$. If none of the candidate colors are $valid$, then by $Fact$ \ref{fact:fact2}, the following Claim is easy to see:

\begin{clm}
\label{clm:clm1}
With respect to the coloring $c_0$, $\forall \beta \in C-F, \exists \alpha \in F_x \cap F_y$ such that there is a $(\alpha,\beta,xy)$ critical path.
\end{clm}

\subsection*{case 2.1:$(S_{xa} \cup S_{xb})\cap F =\emptyset$ }
Since $F=\{1,2,3\}$, $S_{xa} = S_{xb}= \{4,5,6\}$.

\begin{clm}
\label{clm:clm2}
With respect to the coloring $c_0$, all the colors of $C-F$ are weak.
\end{clm}
\begin{proof}
Suppose not. Then there is a strong color in $C-F$. Without loss of generality let $4$ be a strong color. Let $c_0(x,a)=c_0(y,a')=1$ and $c_0(x,b)=c_0(y,b')=2$. Now it is easy to check that the 5-tuple $(x,a,b,{\emptyset},{\emptyset})$ satisfies $configuration\ A$. Let

$$c'_0=ColorExchange(c_0,xa,xb)$$

By $Fact$ \ref{fact:fact4} the color exchange with respect to the edges $xa$ and $xb$ is valid. Thus the coloring $c'_0$ is acyclic.

Since color 4 was strong in coloring $c_0$, there was a $(1,4,xy)$ critical path as well as a $(2,4,xy)$ critical path before $color\ exchange$ (i.e., with respect to the coloring $c_0$). Thus by $Lemma$ \ref{lem:lem2}, $(1,4,xy)$ critical path and $(2,4,xy)$ critical path will not exist after the $color\ exchange$ (i.e., with respect to the coloring $c'_0$). Thus by $Fact$ \ref{fact:fact2}, color 4 is valid for edge $xy$.
\end{proof}

By $Claim$ \ref{clm:clm2}, all the colors of $C-F$ are weak. Each weak color should be actively present in exactly one of $S_{xa}$ or $S_{xb}$. Since there are 3 weak colors, we can infer that either $S_{xa}$ or $S_{xb}$ is such that at least 2 of the weak colors are actively present in it.

\begin{asm}
\label{asm:asm2}
Without loss of generality assume that colors $4$ and $5$ are $actively\ present$ in $S_{xa}$. Let $c(a,k_1) = 4$ and $c(a,k_2) = 5$.
\end{asm}

From $Assumption$ \ref{asm:asm2}, it follows that since $c(x,a) = 1$, there exist $(1,4,xy)$ and $(1,5,xy)$ \emph{critical paths}. The following claim is obvious.

\begin{clm}
\label{clm:clm3}
With respect to the coloring $c_0$, $1 \in S_{ak_1}$ and $1 \in S_{ak_2}$.
\end{clm}

It is easy to verify that the 5-tuple $(x,a,b,{\emptyset},{\emptyset})$ satisfies configuration $A$ with respect to the coloring $c_0$.

$$c_1=ColorExchange(c_0,xa,xb)$$

By $Fact$ \ref{fact:fact4} the color exchange with respect to the edges $xa$ and $xb$ is valid. Thus the coloring $c_1$ is acyclic.

But there were $(1,4,xy)$ and $(1,5,xy)$ \emph{critical paths} before $color\ exchange$ (i.e., with respect to the coloring $c_0$). By $Lemma$ \ref{lem:lem2}, both $(1,4,xy)$ and $(1,5,xy)$ \emph{critical paths} does not exist after the $color\ exchange$ (i.e., with respect to the coloring $c_1$).

Thus even with respect to the coloring $c_1$, if both the colors $4$ and $5$ are not $valid$ for the edge $xy$, by $Fact$ \ref{fact:fact2}, there has to be $(2,4,xy)$ and $(2,5,xy)$ \emph{critical paths}. Thus $2 \in S_{ak_1}$ and $2 \in S_{ak_2}$. Thus we can $Claim$ the following:

\begin{clm}
\label{clm:clm4}
With respect to the coloring $c_1$, $\{1,2\} \subset S_{ak_1}$ and $\{1,2\} \subset S_{ak_2}$. Moreover there will not be any $(1,4,xy)$ and $(1,5,xy)$ \emph{critical paths}.
\end{clm}

Now since the colors 4 and 5 are weak, we try to break the $(2,4,xy)$ and $(2,5,xy)$ \emph{critical paths} by recoloring the edge $xa$.

$$c_2=Recolor(c_1,xa,3)$$

Note that color 3 is a candidate for the edge $xa$ since $S_{xa}= \{4,5,6\}$ and $c(x,b)=1$. And also since $S_{xa} \cap S_{ax} = \emptyset$, by $Lemma$ \ref{lem:lem1} color 3 is $valid$ for the edge $xa$.

Note that with respect to the coloring $c_2$, $F_x \cap F_y = \{1,3\}$. In view of $Claim$ \ref{clm:clm4}, there will not be any $(1,4,xy)$ and $(1,5,xy)$ \emph{critical paths} with respect to the coloring $c_2$ also. If both the colors $4$ and $5$ are not $valid$ for the edge $xy$ still, then by $Fact$ \ref{fact:fact2}, there has to be $(3,4,xy)$ and $(3,5,xy)$ \emph{critical paths} implying $3 \in S_{ak_1}$ and $3 \in S_{ak_2}$. Thus combined with $Claim$ \ref{clm:clm4}, we infer the following:

\begin{clm}
\label{clm:clm5}
With respect to the coloring $c_2$, we have $S_{ak_1} = S_{ak_2} = \{1,2,3\}$. Moreover there will not be any $(1,4,xy)$ and $(1,5,xy)$ \emph{critical paths}.
\end{clm}

Now the 5-tuple $(a,k_1,k_2,\{k_3\},\{x\})$ satisfies configuration $A$.

$$c_3=ColorExchange(c_2,ak_1,ak_2)$$

By fact \ref{fact:fact3} if there is any bichromatic cycle (recalling that $c_3(a,x)=3$), it has to be either a $(5,3)$ or $(6,3)$ bichromatic cycle that passes through vertex $a$ and hence vertex $x$. But any cycle that passes through vertex $x$ should contain edge $xb$ also. Since $c_3(x,b)=1$, this is a contradiction and we infer that $c_3$ is acyclic.

There was a $(3,4,xy)$ critical path as well as a $(3,5,xy)$ critical path before $color\ exchange$ (i.e., with respect to the coloring $c_2$). Thus by $Lemma$ \ref{lem:lem2}, both these critical paths does not exist after the color exchange (i.e., with respect to the coloring $c_3$) (Note that $k_1$, $k_2 \notin \{x,y\}$ since $c(a,k_1)=4$ and $c(a,k_2)=5$ $\notin F_x$ or $F_y$. Therefore we can apply $Lemma$ \ref{lem:lem2})

To summarize, $c_3(x,a)=3$, $c_3(x,b)=1$ and thus $F_x \cap F_y=\{1,3\}$. With respect to the coloring $c_3$, there exist no $(3,4,xy)$ and $(3,5,xy)$ critical paths. Recall that by Claim \ref{clm:clm5}, there won't be any $(1,4,xy)$ and $(1,5,xy)$ critical paths with respect to the coloring $c_2$. It is easy to see that even with respect to the coloring $c_3$, there won't be any $(1,4,xy)$ and $(1,5,xy)$ critical paths.

Thus by $Fact$ \ref{fact:fact2}, color 4 and 5 are valid for edge $xy$.

\subsection*{case 2.2:$(S_{xa} \cup S_{xb})\cap F \neq \emptyset$ }
\begin{asm}
\label{asm:asm3}
Without loss of generality let $S_{xa} \cap F \neq \emptyset$. It follows that one of $\{4,5,6\}$ is missing in $S_{xa}$ since $\vert S_{xa}\vert = 3$. Without loss of generality let it be color 5. Also let $c_0(x,a)=c_0(y,a')=1$ and $c_0(x,b)=c_0(y,b')=2$ and $c_0(y,d')=3$.
\end{asm}

\begin{clm}
\label{clm:clm6}
With respect to the coloring $c_0$, there exists a $(2,5,xy)\ critical\ path$. Thus $5 \in S_{xb}$.
\end{clm}
\begin{proof}
Since color 5 is not valid for the edge $xy$, by $Claim$ \ref{clm:clm1} there has to be a $(1,5,xy)$ critical path or a $(2,5,xy)$ critical path. But by $Assumption$ \ref{asm:asm3}, color $5 \notin S_{xa}$ and hence there can not be a $(1,5,xy)$ critical path. Thus there exists a $(2,5,xy)$ critical path.
\end{proof}

\begin{clm}
\label{clm:clm7}
With respect to the coloring $c_0$, all the colors of $C-F$ are weak.
\end{clm}
\begin{proof}
Suppose not. Then there is at least one strong color in $C-F$. Without loss of generality let 4 be a strong color. Thus we have $4 \in S_{xb}$. Combined with $Claim$ \ref{clm:clm6}, we have:

\begin{equation}
\label{eqn:eqn1}
 \{4,5\} \subset S_{xb}.
\end{equation}

Now let

$$c'_0=Recolor(c_0,xa,5)$$

Note that color $5$ is a candidate for the edge $xa$ since $c_0(x,b)=2$ and $5 \notin S_{xa}$ (by $Assumption$ \ref{asm:asm3}). Now we claim that assigning color 5 to the edge $xa$ can not result in any bichromatic cycle. To see this first note that since any cycle containing the edge $xa$ should also contain the edge $xb$, but $c_0(x,b)=2$ and therefore if a bichromatic cycle gets formed it must be a $(2,5)$ bichromatic cycle, implying that there is a $(2,5,xa)$ critical path. But there is already a $(2,5,xy)$ critical path (by $Claim$ \ref{clm:clm6}) and by $Fact$ \ref{fact:fact1} there can not be a $(2,5,xa)$ critical path, a contradiction. Thus coloring $c'_0$ is acyclic.

Note that with respect to the coloring $c'_0$, color 6 remains to be a candidate color for the edge $xy$. Also note that $F_x \cap F_y = \{2\}$. If the candidate color $6$ is not valid for the edge $xy$, then by $Fact$ \ref{fact:fact2} there has to be a $(2,6,xy)$ critical path and thus $6 \in S_{xb}$. Thus combined with $(\ref{eqn:eqn1})$, we have:

\begin{equation}
\label{eqn:eqn2}
 S_{xb}=\{4,5,6\}
\end{equation}

With respect to the coloring $c_0$, color 4 was strong (assumption) and thus there existed a $(1,4,xy)$ critical path. After recoloring the edge $xa$ with color $5$ (i.e., with respect to the coloring $c_1$), the $(1,4,xy)$ critical path gets curtailed to a $(1,4,y,a)$ maximal bichromatic path without containing the vertex $x$. Moreover note that $(1,4,y,a)$ maximal bichromatic path does not contain the vertex $b$, since if $b$ is in this path, then it is an internal vertex and thus both colors $1,4 \in F_b$, a contradiction ($1 \notin F_b$). Thus we have,

\begin{eqnarray}
\label{eqn:eqn3}
\mbox{\emph{With respect to the coloring $c'_0$, a $(1,4,y,a)$ maximal bichromatic path exists,}} \\ \nonumber
\mbox{\emph{but this path does not contain the vertices $x$ or $b$.}}
\end{eqnarray}

Now with respect to the coloring $c'_0$, $F_x \cap F_y = \{2\}$. Let 

$$c''_0=Recolor(c'_0,xb,1)$$

Note that color 1 is a candidate color for the edge $xb$ since $c'_0(x,a)=5$ and $1 \notin S_{xb}=\{4,5,6\}$. Color 1 is $valid$ for the edge $xb$ because any bichromatic cycle containing edge $xb$ should also contain edge $xa$ and since color $1 \notin S_{xa}$ (Recall that $c_0(x,a)=1$. Thus $1 \notin S_{xa}$ with respect to the coloring $c_0$. Therefore $1 \notin S_{xa}$ with respect to the coloring $c'_0$ also.), such a $(1,5)$ bichromatic cycle can not be formed. Thus $c''_0$ is acyclic. 

Thus with respect to coloring $c''_0$, $F_x \cap F_y = \{1\}$. Now by $(\ref{eqn:eqn3})$, with respect to the coloring $c'_0$, there existed a $(1,4,y,a)$ maximal bichromatic path that does not contain vertex $b$ or $x$. Thus noting that $c''_0$ is obtained just by changing the color of the edge $xb$ to $1$, by $Fact$ \ref{fact:fact1} we infer that $c''_0$ can not contain $(1,4,xy)$ critical path.

Thus by $Fact$ \ref{fact:fact2} color 4 is valid for the edge $xy$.
\end{proof}

\begin{clm}
\label{clm:clm8}
In view of $Assumption$ \ref{asm:asm3}, with respect to the coloring $c_0$, each $\alpha \in \{4,5,6\}$ is $actively\ present$ in $S_{xb}$
\end{clm}
\begin{proof}
Suppose not. By $Claim$ \ref{clm:clm6}, we know that color 5 is $actively\ present$ in $S_{xb}$. Without loss of generality let color 6 be not $actively\ present$ in $S_{xb}$. Therefore color 6 is $actively\ present$ in $S_{xa}$. Now let 

$$c'_0=Recolor(c_0,xa,5)$$

Note that color $5$ is a candidate since $5 \notin \{S_{xa}$ (by $Assumption$ \ref{asm:asm3}) and $c_0(x,b)=2$. Now we claim that assigning color 5 to the edge $xa$ can not result in any bichromatic cycle. To see this first note that since any cycle containing the edge $xa$ should also contain the edge $xb$, but $c_0(x,b)=2$ and therefore if a bichromatic cycle gets formed it must be a $(2,5)$ bichromatic cycle, implying that there is a $(2,5,xa)$ critical path with respect to the coloring $c_0$. But in $c_0$ there is already a $(2,5,xy)$ critical path (by Claim \ref{clm:clm1}) and by $Fact$ \ref{fact:fact1} there can not be a $(2,5,xa)$ critical path, a contradiction. Thus coloring $c'_0$ is acyclic.

Now $F_x \cap F_y = \{2\}$. But in $c_0$, there did not exist a $(2,6,xy)$ critical path since by assumption color 6 is actively present in $S_{xb}$. Thus noting that $c'_0$ is obtained just by changing the color of the edge $xa$ to $5$, we infer that $c'_0$ can not contain $(2,6,xy)$ critical path.

Thus by $Fact$ \ref{fact:fact2} color 6 is valid for the edge $xy$.
\end{proof}

Recall that $c_0(x,b)=c_0(y,b')=2$. In view of $Claim$ \ref{clm:clm8}, with respect to the coloring $c_0$, we have:

\begin{equation}
\label{eqn:eqn4}
S_{xb}=S_{yb'}= \{4,5,6\}
\end{equation}

Let

$$c_1=Recolor(c_0,xb,3)$$

Note that color 3 is a candidate for edge $xb$ since $3 \notin \{S_{xb}=\{4,5,6\}$ (by Claim \ref{clm:clm8}) and $c_0(x,a)=1$. Moreover since $S_{xb} \cap S_{bx} = \emptyset$, by $Lemma$ \ref{lem:lem1} color 3 is also $valid$. Thus the coloring $c_1$ is acyclic.

With respect to the coloring $c_1$, $F_x \cap F_y = \{1,3\}$. In view of $Claim$ \ref{clm:clm7} and $Claim$ \ref{clm:clm8}, $\forall \alpha \in \{4,5,6\}$, $\alpha$ is not $actively\ present$ in $S_{xa}$ and thus $(1,\alpha,xy)$ critical path does not exist with respect to the coloring $c_0$. It is true with respect to the coloring $c_1$ also. Hence if none of the colors from $\{4,5,6\}$ is $valid$ for the edge $xy$ with respect to the coloring $c_1$, then by $Fact$ \ref{fact:fact2} there has to be $(3,4,xy)$, $(3,5,xy)$ and $(3,6,xy)$ \emph{critical paths}. Recalling that by $Assumption$ \ref{asm:asm3} $c(y,d')=3$, we infer that $S_{yd'}=\{4,5,6\}$.

Thus with respect to the coloring $c_1$, we have:

\begin{equation}
\label{eqn:eqn5}
S_{yb'}=S_{yd'}=\{4,5,6\}
\end{equation}

The 5-tuple $(y,b',d',\{a'\},{\emptyset})$ is configuration $A$. Now let

$$c_2=ColorExchange(c_1,yb',yd')$$

By $Fact$ \ref{fact:fact4} the color exchange with respect to the edges $yb'$ and $yd'$ is valid. Thus the coloring $c_1$ is acyclic.

For $\alpha \in \{4,5,6\}$ there was a $(3,\alpha,xy)$ critical path before $color\ exchange$ (with respect to coloring $c_1$). Thus by $lemma$ \ref{lem:lem2}, these critical paths does not exist after the $color\ exchange$ (with respect to coloring $c_2$). Also recall that there was no $(1,\alpha,xy)$ critical path with respect to the coloring $c_1$. Noting that the $color\ exchange$ involved only the colors $2$ and $3$ there is no chance of any $(1,\alpha,xy)$ critical path to get formed with respect to the coloring $c_2$.

Thus by $fact$ \ref{fact:fact2}, color $\alpha$ is valid for edge $xy$.

\subsection*{case 3: $\vert F_x \cap F_y\vert = 1$}

\begin{asm}
\label{asm:asm4}
Without loss of generality let  $F_x = \{1,2\}$ and $F_y = \{1,3,4\}$. Thus $F=\{1,2,3,4\}$. Then $C-F = \{5,6\}$. Let $c_0(x,a)=c_0(y,a')=1$, $c_0(x,b)=2$, $c_0(y,b')=3$ and $c_0(y,d')=4$.
\end{asm}

If none of the colors from $C-F$ are $valid$, then by $Fact$ \ref{fact:fact2}, there exist $(1,5,xy)$ and $(1,6,xy)$ $critical\ paths$. We capture this in the following $claim$:

\begin{clm}
\label{clm:clm9}
With respect to coloring $c_0$, there exist $(1,5,xy)$ and $(1,6,xy)$ $critical\ paths$. Thus $\{5,6\} \subset S_{xa}$ and $\{5,6\} \subset S_{ya'}$.
\end{clm}

\begin{clm}
\label{clm:clm10}
With respect to coloring $c_0$, $\{3,4\}\subset S_{xb}$.
\end{clm}
\begin{proof}
Suppose not. Then at least one of $3$, $4$ is missing in $S_{xb}$. Without loss of generality let $4 \notin S_{xb}$. Recalling that $c_0(x,a)=1$, it follows that color $4$ is a candidate color for the edge $xb$. We claim that there exists a $(1,4,xb$) critical path with respect to the coloring $c_0$. Suppose not. Then let

$$c'_0=Recolor(c_0,xb,4)$$

Clearly $c'_0$ is acyclic since any bichromatic cycle being formed should involve the edge $xa$ as well. But $c'_0(x,a)=1$ and hence a $(1,4)$ bichromatic cycle has to be formed, implying that there is a $(1,4,xb$) critical path, a contradiction to our assumption.

With respect to the coloring $c'_0$, $\vert (F_x \cap F_y)=\{1,4\}\vert = 2$, and by $case$ 2 we will be able to find a valid color for the edge xy.

Thus we can infer that there exists a $(1,4,xb$) critical path with respect to the coloring $c_0$. For a $(1,4,xb$) critical path to exist clearly we should have $4 \in S_{xa}$, since $c_0(x,a)=1$. Combined with $Claim$ \ref{clm:clm9}, we get:

\begin{equation}
\label{eqn:eqn6}
S_{xa}=\{4,5,6\}
\end{equation}

Moreover we have $1 \in  S_{xb}$ with respect to $c_0$ since there is a $(1,4,xb$) critical path. Now let the other two colors in $S_{xb}$ be $\{\alpha,\beta\}$. Then $\gamma \in (\{3,5,6\}-\{\alpha,\beta\})$ is a candidate color for the edge $xb$. Let

$$c''_0=Recolor(c_0,xb,\gamma)$$

We claim that $c''_0$ is acyclic. Otherwise if any bichromatic cycle gets formed with respect to the coloring $c''_0$, then it should be a $(\gamma,1)$ bichromatic cycle since any cycle that contains edge $xb$ should contain edge $xa$ also and $c''_0(x,a)=1$, implying that there exists a $(1,\gamma,xb)$ critical path with respect to the coloring $c_0$. If $\gamma =3$, such a critical path can not exist since $3 \notin S_{xa}$ (by $(\ref{eqn:eqn6})$). On the other hand if $\gamma \in \{5,6\}$, by $Fact$ \ref{fact:fact1}, $(1,\gamma,xb)$ critical path can not exist with respect to the coloring $c_0$ since there is already a $(1,\gamma,xy)$ critical path (by $Claim$ \ref{clm:clm9}). Thus we infer that $c''_0$ is acyclic.

With respect to coloring $c''_0$, if $\gamma=3$, $\vert (F_x \cap F_y)=\{1,3\}\vert = 2$, and by $case$ 2 we will be able to find a valid color for the edge $xy$.

With respect to coloring $c''_0$, if $\gamma \in \{5,6\}$ we have $(F_x \cap F_y) = \{1\}$ and $2 \in C-F$. Thus color 2 is a candidate color for the edge $xy$. Moreover since $S_{xa}=\{4,5,6\}$ (by (6)), there can not be a $(1,2,xy)$ critical path and hence by $Fact$ \ref{fact:fact2}, color 2 is valid for the edge $xy$.
\end{proof}

\begin{clm}
\label{clm:clm11}
With respect to the coloring $c_0$, $S_{xb}= \{3,4,1\}$.
\end{clm}
\begin{proof}
Suppose not. Then in view of Claim \ref{clm:clm10}, we can infer that color $1 \notin S_{xb}$. Recall that by Claim \ref{clm:clm9}, $\{5,6\} \subset S_{xa}$. Let the remaining color in $S_{xa}$ be $\alpha$. Let $\beta \in \{3,4\}-\{\alpha\}$. Now let

$$c'_0=Recolor(c_0,xb,1)$$ and

$$c''_0=Recolor(c'_0,xa,\beta)$$

Note that $c''_0$ is proper since $1 \notin S_{xb}$ (by $Assumption$) and $\beta \notin S_{xa}$, by the definition of $\beta$. The coloring $c''_0$ is acyclic since any cycle containing the edge $xa$ should also contain the edge $xb$ (and vise versa), but $c''_0(x,b)=1$ and therefore if a bichromatic cycle gets formed it must be a $(1,\beta)$ bichromatic cycle, implying that $1 \in S_{xa}$. But this is a contradiction since $1 \notin S_{xa}$ with respect to $c_0$ as $c_0(x,a)=1$ and therefore $1 \notin S_{xa}$ with respect to $c''_0$ also.

Now since $\beta \in \{3,4\}$, we have $\vert (F_x \cap F_y)=\{1,\beta\}\vert = 2$ and thus the situation reduces to $case$ 2, thereby enabling us to find a valid color for the edge $xy$.
\end{proof}

\begin{clm}
\label{clm:clm12}
There is a $(1,2,xy)$ critical path . Thus in combination with $Claim$ \ref{clm:clm9} $S_{xa}= \{5,6,2\}$ , $S_{ya'}= \{5,6,2\}$ with respect to the coloring $c_0$.
\end{clm}
\begin{proof}
Suppose not. Let 

$$c'_0=Recolor(c_0,xb,5)$$

Note that color 5 is a candidate color for the edge $xb$ since, by $Claim$ \ref{clm:clm11}, $S_{xb}=\{3,4,1\}$ and $c_0(x,a)=1$. It is also valid since if there is a bichromatic cycle, then it should contain the edges $xa$ and $xb$ and hence it has to be a $(1,5)$ bichromatic cycle, implying that there exists a $(1,5,xb)$ critical path with respect to the coloring $c_0$. But there can not be a $(1,5,xb)$ critical path (by $Fact$ \ref{fact:fact1}) as there is already a $(1,5,xy)$ critical path (by $Claim$ \ref{clm:clm9}). Thus the coloring $c'_0$ is acyclic.

Now with respect to the coloring $c'_0$, $F_x \cap F_y = \{1\}$. Color 2 is a candidate color for the edge $xy$ since $2 \notin (F_x \cup F_y =\{1,3,4,5\})$. Since there is no $(1,2,xy)$ critical path (by assumption), by $Fact$ \ref{fact:fact2}, color 2 is valid for the edge $xy$.
\end{proof}

Recall that $N_{G'}(a)=\{x,k_1,k_2,k_3\}$ and $N_{G'}(b)=\{x,l_1,l_2,l_3\}$. Also recall that by Assumption \ref{asm:asm4}, $c_0(x,a)=c_0(y,a')=1$,$c_0(x,b)=2$,$c_0(y,b')=3$ and $c_0(y,d')=4$. By Claim \ref{clm:clm11} and Claim \ref{clm:clm12}, $S_{xa}= \{5,6,2\}$ and $S_{xb}= \{3,4,1\}$. We make the following $Assumption$:

\begin{asm}
\label{asm:asm5}
Without loss of generality let $c_0(a,k_1)= 5$, $c_0(a,k_2)= 6$, $c_0(a,k_3)= 2$, $c_0(b,l_1)= 3$, $c_0(b,l_2)= 4$ and $c_0(b,l_3)= 1$.
\end{asm}

The main intention of the next two $Claims$ is to establish that $S_{bl_1} = S_{bl_2} =\{2,5,6\}$.

\begin{clm}
\label{clm:clm13}
With respect to the coloring $c_0$, there exist $(2,3,xa)$ and $(2,4,xa)$ critical paths. Thus $ 2 \in S_{bl_1}$, $ 2 \in S_{bl_2}$.
\end{clm}
\begin{proof}
Suppose not. Then without loss of generality let there be no $(2,3,xa)$ critical path. Let

$$c'_0=Recolor(c_0,xa,3)$$

Note that color 3 is a candidate color for edge $xa$ since $3 \notin (S_{xa} =\{2,5,6\})$ (by Claim \ref{clm:clm12}) and $c_0(x,b)=2$. It is also valid since if there is any bichromatic cycle containing edge $xa$, then it should also contain edge $xb$ and since $c_0(x,b)=2$, it has to be a $(2,3)$ bichromatic cycle, implying that there is a $(2,3,xa)$ critical path, a contradiction to our assumption. Thus the coloring $c'_0$ is acyclic.

With respect to the coloring $c'_0$, $c'_0(y,b')=3$ and $(F_x \cap F_y) = \{3\}$. Now if one of the colors 5 and 6 are valid for the edge $xy$, we are done. Otherwise by $Fact$ \ref{fact:fact2}, there are $(3,5,xy)$ and $(3,6,xy)$ critical paths. Thus

\begin{equation}
\label{eqn:eqn7}
\{5,6\} \subset S_{yb'}
\end{equation}

Let,

$$c''_0=Recolor(c'_0,xb,5)$$

First note that color 5 is a candidate for the edge $xb$ since $5 \notin (S_{xb}=\{3,4,1\})$ (by $Claim$ \ref{clm:clm11}) and $c'_0(x,a)=3$ . It is also valid since if there is any bichromatic cycle containing the edge $xb$ then it should also contain edge $xa$ and since $c'_0(x,a)=3$, it has to be a $(3,5)$ bichromatic cycle,implying that there exists a $(3,5,xb)$ critical path. But there can not be a $(3,5,xb)$ critical path (by $Fact$ \ref{fact:fact1}) as there is already a $(3,5,xy)$ critical path. Thus the coloring $c''_0$ acyclic.

Now with respect to the coloring $c''_0$, $(F_x \cap F_y) = \{3\}$ and $2 \notin (F_x \cup F_y)=\{1,3,4,5\}$. Color 2 is a $candidate$ for the edge $xy$. If it is $valid$ then we are done. Otherwise by $Fact$ \ref{fact:fact2}, there exists a $(3,2,xy)$ critical path.

Thus $2 \in S_{yb'}$ and in combination with $(\ref{eqn:eqn7})$, we get,

\begin{equation}
\label{eqn:eqn8}
S_{yb'}=\{2,5,6\}
\end{equation}

Recall that $S_{ya'}=\{2,5,6\}$ by $Claim$ \ref{clm:clm12} with respect to the coloring $c_0$. It is easy to see that $S_{ya'}=\{2,5,6\}$ even with respect to the coloring $c_2$. Now in view of Assumption \ref{asm:asm4}, we have the 5-tuple $(y,a',b',\{d'\},{\emptyset})$ in $Configuration\ A$. Let,

$$c'''_0=ColorExchange(c''_0,ya',yb')$$

By $Fact$ \ref{fact:fact4}, the color exchange with respect to the edges $ya'$ and $yb'$ is valid. Thus the coloring $c'''_0$ is acyclic.

There was a $(3,6,xy)$ critical path before $color\ exchange$ (i.e., with respect to the coloring $c''_0$) since otherwise color 6 would have been valid for the edge $xy$ with respect to the coloring $c''_0$. Thus by $Lemma$ \ref{lem:lem2} no $(3,6,xy)$ critical path exists after the $color\ exchange$ (i.e., with respect to the coloring $c'''_0$). Thus by $Fact$ \ref{fact:fact2}, color $6$ is valid for edge $xy$.
\end{proof}

\begin{clm}
\label{clm:clm14}
With respect to the coloring $c''_0$, $\forall \alpha \in \{3,4\}$ and $\forall \beta \in \{5,6\}$, there exist $(\alpha,\beta,b,a)$ maximal bichromatic path which ends at vertex $a$ with an edge colored $\beta$. Thus $S_{bl_1}=\{2,5,6\}$ and $S_{bl_2}=\{2,5,6\}$.
\end{clm}
\begin{proof}
Suppose not. Then $\exists \alpha \in \{3,4\}$ and $\exists \beta \in \{5,6\}$ such that there is no $(\alpha,\beta,b,a)$ maximal bichromatic path which ends at vertex $a$ with an edge colored $\beta$. Without loss of generality let $\alpha = 3$ and $\beta = 5$. Now let,

$$c'_0=Recolor(c_0,xa,3)$$ and

$$c''_0=Recolor(c'_0,xb,5)$$

Note that $c''_0$ is a proper coloring ( since ($3 \notin S_{xa}=\{2,5,6\}$ and $c''_0(x,b)=5$) and ($5 \notin S_{xb}=\{3,4,1\}$ and $c''_0(x,b)=3$ )). Now to see that $c''_0$ is acyclic, note that if there is a bichromatic cycle with respect to the coloring $c''_0$, then it should contain both the edges $xa$ and $xb$, thus forming $(3,5)$ bichromatic cycle, implying that there should be a $(3,5,a,b)$ maximal bichromatic path which ends at vertex $a$ with an edge colored $3$ with respect to the coloring $c_0$, a contradiction to our assumption.

Note that with respect to the coloring $c''_0$, $F=\{1,3,4,5\}$ and thus color 2 is a candidate color for the edge $xy$. By $Claim$ \ref{clm:clm13} there was a $(2,3,xa)$ critical path with respect to the coloring $c_0$. From this it is easy to see that with respect to the coloring $c''_0$, there is a $(3,2,xb)$ critical path. Thus by $Fact$ \ref{fact:fact1} there can not be a $(3,2,xy)$ critical path with respect to the coloring $c''_0$. Hence color 2 is valid for the edge $xy$.

Thus $\forall \alpha \in \{3,4\}$ and $\forall \beta \in \{5,6\}$, there exist $(\alpha,\beta,b,a)$ maximal bichromatic path which ends at vertex $a$ with an edge colored $\beta$. Thus recalling that $c_0(b,l_1)=3$ and $c_0(b,l_2)=4$ with respect to the coloring $c_0$, we have,

\begin{eqnarray}
{\label{eqn:eqn9} \{5,6\} \subset S_{bl_1} }\\
{\label{eqn:eqn10} \{5,6\} \subset S_{bl_2} }
\end{eqnarray}

By $Claim$ \ref{clm:clm13}, $2 \in S_{bl_1}$ and $2 \in S_{bl_2}$. Thus we have,

\begin{equation}
\label{eqn:eqn11}
S_{bl_1}=S_{bl_2}=\{2,5,6\}
\end{equation}

\end{proof}

Now let,

$$c_1=Recolor(c_0,xb,5)$$

Recalling Claim \ref{clm:clm11}, $s_{xb}=\{3,4,1\}$ and $c_0(x,a)=1$, color 5 is a candidate for the edge $xb$. Moreover color 5 is also valid since if there is any bichromatic cycle containing the edge $xb$ then it should also contain edge $xa$ and since $c_0(x,a)=1$, it has to be a $(1,5)$ bichromatic cycle,implying that there exists a $(1,5,xb)$ critical path with respect to the coloring $c_0$. But there can not be a $(1,5,xb)$ critical path (by $Fact$ \ref{fact:fact1}) as there is already a $(1,5,xy)$ critical path (by $Claim$ \ref{clm:clm9}). Thus the coloring $c_1$ is acyclic.

Recall that by $Claim$ \ref{clm:clm14}, with respect to the coloring $c_0$, there was a $(3,5,b,a)$ maximal bichromatic path that ends at vertex $a$ with an edge colored $5$. After the recoloring of edge $xb$ with color $5$ (i.e., with respect to the coloring $c_1$),it is easy to see that this $(3,5,b,a)$ maximal bichromatic path gets extended to a $(3,5,xa)$ critical path. Thus we have,

\begin{equation}
\label{eqn:eqn12}
\mbox{\emph{With respect to the coloring $c_1$, there exists a $(3,5,xa)$ critical path.}}
\end{equation}

Recall that by $Claim$ \ref{clm:clm13}, with respect to the coloring $c_0$, there existed a $(2,3,xa)$ critical path. After recoloring the edge $xb$ with color $5$ (i.e., with respect to the coloring $c_1$), the $(2,3,xa)$ critical path gets curtailed to a $(2,3,a,b)$ maximal bichromatic path that ends at vertex $b$ with an edge colored $3$. Note that $(2,3,a,b)$ maximal bichromatic path does not contain the vertex $y$, since if $y$ is in this path, then it is an internal vertex and thus both colors $2,3 \in F_y$, a contradiction ($2 \notin F_b$). Thus noting that $c_1(b,l_1)=3$, we have,

\begin{eqnarray}
\label{eqn:eqn13}
\mbox{\emph{With respect to the coloring $c_1$, there exists a $(2,3,a,b)$ maximal bichromatic path that ends at vertex $b$}} \\ \nonumber
\mbox{\emph{with an edge colored 3. This path contains the edge $bl_1$ but does not contain vertex $y$.}}
\end{eqnarray}

In view of $Claim$ \ref{clm:clm14}, we have $S_{bl_1} = S_{bl_2} = \{2,5,6\}$. The 5-tuple $(b,l_1,l_2,\{l_3\},\{x\})$ is in $configuration\ A$. Let,

$$c_2=ColorExchange(c_1,bl_1,bl_2)$$

By $Fact$ \ref{fact:fact3} if there is any bichromatic cycle, recalling that $c_2(x,b)=5$, there has to be either $(3,5)$ or $(4,5)$ bichromatic cycle that passes through vertex $x$. But any cycle that passes through vertex $x$ should contain edge $xa$ also. Since $c_2(x,a)=1$, this is a contradiction and we infer that $c_2$ is acyclic.

Note that by $(\ref{eqn:eqn13})$ there existed  $(2,3,a,b)$ maximal bichromatic path containing the edge $bl_1$ with respect to the coloring $c_1$. Since the color of edge $bl_1$ is changed in $c_2$, this path gets curtailed to a $(2,3,a,l_1)$ maximal bichromatic path which now ends at the vertex $l_1$ since $3 \notin F_{l_1}$ with respect to the coloring $c_2$. Note that it still does not contain vertex $y$. Thus we have,

\begin{equation}
\label{eqn:eqn14}
\mbox{\emph{With respect to the coloring $c_2$, there exists a $(2,3,a,l_1)$ maximal bichromatic path which does not contain vertex $y$.}}
\end{equation}

But before $color\ exchange$ (i.e., with respect to the coloring $c_1$) by $(\ref{eqn:eqn12})$ there was a $(3,5,xa)$ critical path. Clearly this path passes through the vertex $b$. Thus by $Lemma$ \ref{lem:lem2}, the $(3,5,xa)$ critical path, does not exist after the color exchange (with respect to the coloring $c_2$) (It easy to see that $l_1$, $l_2 \notin \{x,a\}$ since $1 \notin F_{l_1}$, $F_{l_2}$ but $1 \in F_x$, $F_a$. Therefore $Lemma$ \ref{lem:lem2} can be applied). Thus we have,

\begin{equation}
\label{eqn:eqn15}
\mbox{\emph{With respect to the coloring $c_2$, there does not exists any $(3,5,xa)$ critical path.}}
\end{equation}

Now let

$$c_3=Recolor(c_2,xa,3)$$

By $Claim$ \ref{clm:clm12}, $s_{xa}=\{2,5,6\}$ with respect to the coloring $c_0$ and $s_{xa}=\{2,5,6\}$ even with respect to the coloring $c_2$. Thus color 3 is candidate for edge $xa$ since $3 \notin S_{xa}$ and $c_2(x,b)=5$. Coloring $c_3$ is also acyclic since if there is any bichromatic cycle containing edge $xa$ then it should also contain edge $xb$. But $c_3(x,b)=5$ and $c_3(x,a)=3$. Thus it has to be a $(3,5)$ bichromatic cycle, implying that there exists a $(3,5,xa)$ critical path with respect to the coloring $c_2$, a contradiction (by $(\ref{eqn:eqn15})$).

Note that by $(\ref{eqn:eqn14})$ there existed  $(2,3,a,l_1)$ maximal bichromatic path with respect to the coloring $c_2$. Since the color of edge $xa$ is changed in $c_3$ to color 3, it is easy to see that this path gets extended to a $(3,2,x,l_1)$ maximal bichromatic path which now starts at the vertex $x$ since $2 \notin F_x$ with respect to the coloring $c_3$. Note that it still does not contain vertex $y$.

Now with respect to the coloring $c_3$, $F = \{1,3,4,5\}$ and $F_x \cap F_y = \{3\}$. Thus color 2 is a candidate for the edge $xy$. Since $(2,3,x,l_1)$ maximal bichromatic path contains vertex $x$ and does not contain vertex $y$, by $Fact$ \ref{fact:fact1} there can not be $(2,3,xy)$ critical path. Thus by $Fact$ \ref{fact:fact2} color 2 is valid for the edge $xy$.

\end{proof}


\end{document}